\documentclass[a4wide]{amsart}
  \usepackage[utf8]{inputenc}
  \usepackage{amssymb,exscale,bbm,xpatch,mathrsfs,nicefrac}
  \usepackage[hidelinks]{hyperref}
  \usepackage{enumitem}
  \usepackage{tikz} \usepackage{tikz-cd} \usepackage{graphicx}

 \newtheorem{theorem}{Theorem}[section]
 \newtheorem{coro}[theorem]{Corollary}
 \newtheorem{lemma}[theorem]{Lemma}
 \newtheorem{proposition}[theorem]{Proposition}
 \theoremstyle{definition}
 \newtheorem{definition}[theorem]{Definition}
 \newtheorem{remark}[theorem]{Remark}

 \newtheorem{example}[theorem]{Example}
 
 \numberwithin{equation}{section}

\newcommand{\norm}[1]{\| #1 \|}
\newcommand{\bignorm}[1]{\bigl\| #1 \bigr\|}
\newcommand{\Bignorm}[1]{\Bigl\| #1 \Bigr\|}

\newcommand {\sfrac}[2] { {{}^{#1}\!\!/\!{}_{#2}}}

\newcommand{\eins}{\mathbbm 1}

\DeclareMathOperator{\calO}{\mathcal{O}}

\newcommand{\NN}{\mathbb{N}}

\newcommand{\RR}{\mathbb{R}}

\newcommand{\TT}{\mathbb{T}}
\newcommand{\DD}{\mathbb{D}}

\newcommand{\vphi}{\varphi}

\newcommand{\Sum}[2][\relax]{%
 \ifx#1\relax \sideset{}{_{#2}}\sum 
 \else \sideset{}{^{#1}_{#2}}\sum
 \fi}

\newcommand{\EE}{\mathbb E}
\newcommand{\PP}{\mathbb P}

\newcommand{\CM}{\kappa}
\newcommand{\MF}{\mathcal M}
\newcommand{\CS}{\mathcal C}
\newcommand{\BOUNDED}{\mathscr L}

\newcommand{\om}{\omega}
\newcommand{\Om}{\Omega}

\newcommand{\BL}{\mathcal{L}}

\makeatletter
\newcommand{\sprod}[2]{\ensuremath{%
  \setbox0=\hbox{\ensuremath{#2}}
  \dimen@\ht0
  \advance\dimen@ by \dp0
  \left[ #1\rule[-\dp0]{0pt}{\dimen@}, #2\hspace{1pt}\right]}}
\newcommand{\bsprod}[2]{\ensuremath{%
  \setbox0=\hbox{\ensuremath{#2}}
  \dimen@\ht0
  \advance\dimen@ by \dp0
  \bigl[ #1\rule[-\dp0]{0pt}{\dimen@}, #2\hspace{1pt}\bigr]}}

 \makeatother
\newcommand{\dprod}[2]{\ensuremath{\langle#1,#2\rangle}}

\newcommand{\SUCHTHAT}{:\; }

\newcounter{aufzi}

\newcounter{aufzii}

\newcounter{aufziii}

\newcommand{\dsp}{\displaystyle}
\newcommand{\proofskip}{\medskip\noindent}

\begin{document}

\title[Representing families in Hardy spaces]{Representing families of reproducing kernels in Hardy spaces}

\author[Bernhard H. Haak]{Bernhard H. Haak}
\address{%
  Institut de Math\'ematiques de Bordeaux\\
  Universit\'e de Bordeaux\\
  351, cours de la Lib\'eration\\
  33405 Talence cedex\\
  France}
\email{bernhard.haak@math.u-bordeaux.fr}

\thanks{\today}

\subjclass{Primary 30H10;  Secondary 42C15, 30E10, 46E30}

\keywords{Hardy spaces, reproducing kernels,  representing systems, non-tangential sampling,
  mixed-norm coefficient spaces,  Fejér-type approximation}

\date{\today}


\begin{abstract}
  We study representing systems of Cauchy kernels in the Hardy spaces
  $H^p(\DD)$, $1<p<\infty$. Using an abstract quotient-space
  principle and non-tangential maximal functions, we introduce a
  coefficient space based on the shadows of kernel points and
  characterize when its synthesis operator is surjective. In contrast
  to earlier constructions based on radial blockings, our approach
  allows representing systems with vanishing upper density, answering
  a question raised by Baranov and Batenev. We also investigate
  representations arising from randomly rotated configurations and
  Poisson point processes. Finally, we construct uniformly bounded
  linear sampling operators and show that, after passing to a single
  inductively selected subsequence, their Cesàro averages converge
  strongly to the identity on every $H^p$, $1<p<\infty$. The
  subsequence is independent of both the function and the exponent
  $p$.
\end{abstract}

\maketitle

\allowdisplaybreaks

\section{Introduction}\label{intro}

Throughout the paper we denote by $H^p=H^p(\DD)$ the Hardy space on
the unit disc, and by $m$ the normalized arc measure on
$\TT=\partial\DD$.  For $\lambda\in\DD$, put
\[
  k_\lambda(z):=\frac{1}{1-\overline\lambda\, z},
\]
the Cauchy (or Szeg\H{o}) kernel. It is well known (see e.g.
\cite{BaranovBatenev:2023}) that normalized Cauchy kernels cannot form
a frame for $H^2$. Indeed, the upper (Bessel) inequality
\[
  \sum_{\lambda \in \Lambda}  \frac{1}{\norm{k_\lambda}^2} \bigl| [ f, k_\lambda ] \bigr|^2 \le B \norm{ f }_{H^2}^2,
\]
tested on the constant function $f=\eins$, would imply $\Lambda$ to be
Blaschke, so that a nonzero Blaschke product vanishes at every point
of $\Lambda$: hence the lower inequality fails. This excludes Riesz
bases as well, although suitable subfamilies can form Riesz sequences.

\medskip

Nevertheless, the kernels can form \emph{representing systems}. By
this we mean that, for a fixed enumeration of a sequence $\Lambda
\subset\DD$, every $f\in H^p$ admits scalar coefficients
$(c_\lambda)$ such that
\begin{equation}\label{eq:representations}
  f = Sc := \sum_{\lambda\in\Lambda}c_\lambda k_\lambda,
\end{equation}
where the series converges in the $H^p$-norm, without requiring
uniqueness of the coefficients. The existence question for $H^2$ was
raised by Fricain--Khoi--Lef\`evre \cite{FricainKhoiLefevre} and
answered by Speransky--Terekhin \cite{SperanskyTerekhin:2018}, using
discrete Fourier transform in $H^2$, and extended to radial systems in
\cite{SperanskyTerekhin:2019} using frame theory. Later,
Baranov--Batenev \cite{BaranovBatenev:2023} extended the picture to
$H^p$, $1<p<\infty$, giving both an elementary proof based on
discretization of the Cauchy formula and a frame-theoretic argument
for more general radial configurations.

\bigskip

\noindent There is a simple dual interpretation of \eqref{eq:representations}.
If $1<p<\infty$, $q$ is its conjugate exponent, and $g\in H^q$, then,
with the usual Hardy-space duality pairing,
\[
  \langle \ell_g,Sc\rangle
    =\sum_{\lambda\in\Lambda}c_\lambda\overline{g(\lambda)}.
\]
Thus synthesis by kernels is naturally related to quadrature
functionals on $H^q$. The cited constructions
\cite{SperanskyTerekhin:2018,SperanskyTerekhin:2019,
  BaranovBatenev:2023} approximate the Cauchy integral by sampling on
circles of radius $r\uparrow1$. They also suggest coefficient spaces
organized in radial blocks. Since coefficients need not be unique, the
appropriate norm for a representation is a quotient norm, which
motivates the abstract principle developed in Section~\ref{sec:FA}
below.

\medskip

\noindent One might expect the radial quadrature viewpoint to impose a minimum
density near the boundary. Baranov--Batenev explicitly asked whether a
representing system could have vanishing upper density, and suggested
a negative answer. We shall give a positive answer. Our approach uses
a second description of Hardy spaces, through non-tangential boundary
limits and maximal functions, in place of a fixed radial
organization. The resulting \emph{shadow coefficient space}
$\CS_\Lambda$ provides a geometric characterization of bounded,
unconditionally convergent kernel representations
(Theorem~\ref{thm:nt-representation}). It is striking that this leads
to a different coefficient space
\[
\underbrace{\ell^1\bigl(\NN;\ell^p(\Lambda_n)\bigr)}_{\text{radial blockings}}
\qquad\text{---versus---}\qquad
\underbrace{L^p\bigl(\TT;\ell^1(\Lambda \cap \Gamma(\zeta))\bigr)}
_{\text{non-tangential geometry}}.
\]
(where $\Gamma(\zeta)$ denotes a Stolz region with vertex
$\zeta \in \TT$), thus in an exchange of the order of the $\ell^1$ and
$L^p$ norms. As mentioned above, our approach does allow to treat
sequences of zero upper density naturally, and also gives
representations for randomly rotated configurations and rotationally
invariant Poisson point processes.

\medskip

\noindent Finally, Section~\ref{sec:explicit} returns to the quadrature
interpretation and asks how explicitly the coefficients can be
chosen. We construct linear sampling coefficient maps
$R_n:H^p\to\CS_\Lambda$, given by the same pointwise formulas
independently of $1<p<\infty$, such that
\[
  S R_n f\rightharpoonup f,\qquad   \sup_n\norm{R_n}_{\BOUNDED(H^p,\CS_\Lambda)}<\infty.
\]
For equally spaced radial configurations, ordinary Ces\`aro averages
converge strongly (Theorem~\ref{thm:fejer-std-config}). In the general
case Theorem~\ref{thm:Fejer-on-subsequences} shows that we can
inductively select a subsequence of $(R_n)$ which has the property
that
\[
  \frac1N\sum_{j=1}^N SR_{\vphi(j)}f\longrightarrow f   \qquad\text{in }H^p,\quad 1<p<\infty.
\]
The result gives universal linear \emph{approximations} (but not exact
representations) from sampled values. So the same sequence of finite
sampling operators provides strong approximation simultaneously for
every function in every $H^p$, $1<p<\infty$, without any
function-dependent selection or abstract convexification.  Strong
reconstruction from sampled values is of course not achieved for the
first time, see for example \cite{CimaStressin}.  Our approximants
come from geometric shadow quadrature, with uniform coefficient
control and a subsequence chosen independently of both $f$ and $p$,
followed by ordinary Cesàro averaging.

\section{Functional analytic framework}\label{sec:FA}
The existing radial constructions and our geometric approach share a
simple functional-analytic principle: embeddings of dual spaces
produce quotient representations.

\begin{lemma}[Abstract  quotient principle]\label{lem:abstract-realisation}
  Let $X$ and $Z$ be Banach spaces and let $T : X' \to Z'$
be an isomorphic embedding, i.e. there exist constants $c,C>0$
such that
\[
  c \, \norm{x'}_{X'} \le \norm{Tx'}_{Z'} \le C \, \norm{x'}_{X'}
\]
for all $x' \in X'$. Assume moreover $T'(Z) \subset X$, where we
identify $X$ and $Z$ with their canonical images in their respective
biduals. Then $S := T'|_{Z} : Z \to X$ is a bounded surjection and setting 
$N = \ker(S)$, we have $Z/N \simeq X$ with
\[
  \frac{1}{C} \, \norm{x}_X \; \le \; \inf\bigl\{  \norm{z}_Z \SUCHTHAT z \in Z, Sz = x \bigr\} \; \le \; \frac{1}{c}\, \norm{x}_X
\]
\end{lemma}
\begin{proof}
  By assumption $S \in \BL(Z,X)$ and $S'=T$. Since $T$ is bounded
  below, $S'$ has closed range and trivial kernel. The closed range
  theorem therefore implies that $S$ has closed and dense range, hence
  is surjective. The norm estimates are now standard: the adjoint of
  the induced isomorphism $\overline{S}: Z/N \to X$ identifies with
  $T$. The lower bound for $T$ then implies
  $\norm{ (\overline{S}')^{-1} } \le \frac{1}{c}$ and the claimed
  upper bound follows by exchanging adjoint and inverse. Finally, if
  $Sz = x$, then
  $|\dprod{x}{x'}| = |\dprod{z}{T x'}| \le \norm{z} C
  \norm{x'}$. Taking supremum over the unit ball of $X'$ shows the
  claimed lower bound.
\end{proof}

We will now specialize the lemma to a setup of Banach function
spaces. This approach is vaguely inspired by a result of Bonsall
\cite{Bonsall:ell1-representations}, which is concerned with the case
$Z=\ell_1$ only. Let $(\Om, \Sigma, \mu)$ be a $\sigma$-finite measure
space. The space of (equivalence classes of) $\mu$-measurable
functions over $\Om$ is denoted by $L^0(\Om)$. A Banach function space
is a subspace $Z\subseteq L^0(\Om)$ such that
  \begin{enumerate}[label=(\roman*)]
  \item\label{BFS-item-i}
    If $f\in Z$ and $g\in L^0(\Omega)$ with $|g|\le|f|$ a.e.,
    then $g\in Z$ and $\norm{ g }_{Z}\le\norm{ f }_{Z}$.
  \item\label{BFS-item-ii}
    For every measurable $E\subseteq\Omega$ of positive measure,
    there exists a measurable subset $F\subseteq E$ of positive
    measure with $\eins_F\in Z$.
\end{enumerate}
We will require our Banach function spaces to have absolutely
continuous norms, that is
\begin{enumerate}[label=(\roman*),resume]
\item\label{BFS-item-iii}
    If $f\in Z$ and $(E_n)$ is a decreasing sequence of measurable
    subsets with $\eins_{E_n}\to0$ a.e., then
    \(\norm{ f\,\eins_{E_n} }_{Z}\to0\).
\end{enumerate}
By \cite[Proposition~3.6, Theorem~3.11]{Bennett-Sharpley}, in Banach
function spaces with absolutely continuous norms a suitable version of
dominated convergence holds, and simple functions are dense.  The {\it
  Köthe dual} $Z' \subset L_0(\Omega)$ of a Banach function space $Z$
is defined as
\[
  Z' = \Bigl\{ g \in L^0(\Omega) \SUCHTHAT
   \norm{ g }_{Z'} = \sup\left\{\left|\int_\Omega g\, f \,d\mu\right| \SUCHTHAT \norm{ f }_Z \leq 1\right\} < \infty \Bigr\}
\]
Since $Z$ has absolutely continuous norm, its Banach dual
is canonically and isometrically identified with its Köthe dual.

\medskip

\begin{theorem}[Representations over Banach function spaces]\label{thm:bonsall}
  Let $Z$ be a Banach function space with absolutely continuous norm,
  $X$ a Banach space and \( \vphi: \Om \to X \) a strongly measurable
  function. Assume that
  \[
    \langle x', \vphi(\cdot) \rangle  \in Z'
  \]
  for all $x' \in X'$, and that for some constants $M_1, M_2>0$,   
  \begin{enumerate}[label=(\alph*), ref=(\alph*),leftmargin=*]
\item\label{item:Bonsall-Z-upper} For all $x' \in X'$, 
  $\bignorm{  \dprod{x'}{\vphi(\cdot)} }_{Z'} \le M_1 \norm{x'}_{X'} $
\item\label{item:Bonsall-Z-lower} For all $x' \in X'$, 
  \( \bignorm{ \dprod{x'}{\vphi(\cdot)} }_{Z'} \ge M_2 \norm{x'}_{X'} \).
\end{enumerate}
Then any $x \in X$ can be written as a Pettis integral
\begin{equation} \label{eq:repres}
  x = \int_\Om  z(t) \, \vphi(t)\,d\mu(t)
\end{equation}
with coefficient function $z \in Z$. Moreover,
\begin{equation} \label{eq:norm-estimate} \frac{1}{M_1} \norm{x}_X \le
  \; \inf \; \left\{ \norm{ z }_Z \SUCHTHAT \quad x = \int_\Om z(t) \,
    \vphi(t)\,d\mu(t) \right\}\; \le \frac{1}{M_2} \norm{x}_X
\end{equation}
\end{theorem}
\begin{proof}
  Let us define $T: X' \to Z'$ by $T x' =
  \dprod{x'}{\vphi(\cdot)}$. By our hypotheses, $T$ is an isomorphic
  embedding, so we have to check that $T'(Z) \subseteq X$. To this
  end, fix $z\in Z$ and choose increasing measurable sets $F_n \uparrow \Omega$
  of  finite measure. Then set
  \[
    E_n = F_n \cap \{ |z| \le n \} \cap \{ \norm{ \vphi }_X \le n \}
  \]
  Then $E_n$ is measurable, the sequence $(E_n)$ is increasing, and
  $\eins_{E_n} \uparrow \eins$ almost everywhere.  The absolute
  continuity of the norm of $Z$ gives $z_n := z \eins_{E_n} \to
  z$. Further, $\norm{z(t) \vphi(t) }_X \le n^2$ for $t \in E_n$, so
  that $z_n(\cdot) \vphi(\cdot) \in L^1(\Omega; X)$. Let
  \[
    x_n := \int_{\Omega} z_n(t)\,\vphi(t)\,d\mu(t)
  \]
  Then $\dprod{x'}{x_n}_{X' \times X} = \dprod{Tx'}{z_n}_{Z' \times Z}$ and so 
  \[
    \norm{ x_n - x_m }_X = \sup_{ \norm{x'} \le 1} \bigl| \dprod{ Tx'}{z_n - z_m}\bigr|  \le M_1 \norm{ z_n - z_m }_Z,
  \]
  showing that $(x_n)$ converges to some  $x \in X$. Moreover
\[
    \Bigl|\dprod{ Tx'}{z_n - z }\Bigr| \; \le  \; \norm{ z_n - z }_Z \, \norm{ T x'}_{Z'}
    \le  M_1 \norm{z_n - z}_Z \, \norm{x'}_{X'} \to 0    
\]
  and consequently
  \[
    \dprod{x'}{x} = \lim_n \;\dprod{x'}{x_n}
                  = \lim_n \; \dprod{ T x'}{z_n}
                  = \dprod{T x'}{z} 
  \]
  so that $T' z = x \in X$ as desired.  Lemma~\ref{lem:abstract-realisation}
  applies and yields $S: Z \to X$. By construction
  \[
    S (\eins_A z) = \int_A z(t) \vphi(t)\,d\mu(t)
  \]
  for all measurable $A \subset \Omega$, so \eqref{eq:repres} holds as a
  Pettis integral and \eqref{eq:norm-estimate} follows.
\end{proof}

In Section~\ref{ss:quotient-example} we illustrate why the quotient
structure is necessary rather than merely a matter of notation.

\section{Deterministic representations in Hardy spaces}

The rest of this paper is concernend with $H^p$ spaces on the disc, we
refer the reader to \cite{Duren:Hp-spaces} as one standard reference
for the Hardy-space background, boundary limits, duality and
non-tangential maximal estimates.

\subsection{A classical radial example}\label{ss:ST-example}
As a first illustration of Theorem~\ref{thm:bonsall}, we recover an
$H^p$ version of a radial construction of Speransky--Terekhin
\cite{SperanskyTerekhin:2018,SperanskyTerekhin:2019}, also covered in
greater generality by Baranov--Batenev \cite{BaranovBatenev:2023} (see
Theorem~\ref{thm:BB} below).  Recall the reproducing kernel function
$k_\lambda(z) = \frac{1}{1-\overline{\lambda}\, z}$.  Let
$\lambda_{n,k}= r_n\om_n^k$, where $r_n=1-\sfrac{1}{n}$,
$\om_n=e^{2\pi i/n}$, and
  \[
    \Lambda = \{ \lambda_{k,n} \SUCHTHAT n\ge 2, k=0\ldots{}n{-}1 \}.
  \]
  Then any $f \in X = H^p$ can be written as
\[
  f(z) = \sum_{\lambda \in \Lambda} z_\lambda \, \frac{ k_\lambda(z) }{\norm{k_\lambda}_{H^p}}
\]
where $(z_\lambda)$ satisfies
\[
  \sum_{n=2}^\infty \Bigl( \sum_{|\lambda|=r_n} |z_{\lambda}|^p  \Bigr)^{\sfrac1p}
  \simeq
  \norm{ f }_{H^p}
\]
with uniform constants.  

\begin{proof}
  Let $Z$ be the space of scalar sequences over $\Lambda$ such that
  \[
    \norm{z} = \sum_{n=2}^\infty  \Bigl(  \sum_{k=0}^{n-1}    |z_{n,k}|^p \Bigr)^{\sfrac1p} 
  \]
  is finite. With respect to the usual dual pairing, $Z'$ consists of
  scalar sequences for which 
  \[
   \dsp \norm{w}_{Z'} = \sup \Bigl\{   \Bigl(   \sum_{k=0}^{n-1} |w_{n,k}|^q \Bigr)^{\sfrac1q} \SUCHTHAT \quad n\ge 2 \Bigr\}
 \]
 is finite, where $q$ is the dual index to $p$, i.e. $\frac1p+\frac1q = 1$.  

 \proofskip Put $\vphi(\lambda):=k_\lambda/\norm{k_\lambda}_{H^p}$.
 Now let $g\in H^q$, identified with an element of $X'$. Recall that
 $\norm{ k_{\lambda_{n,k}} }_{H^p} \simeq n^{\sfrac1q}$ by Lemma
 \ref{lem:Hp-norms}. Therefore,
 \[
    \sum_{k=0}^{n-1} |\dprod{g}{\vphi(\lambda_{n,k})}|^q
   \asymp \frac1n\sum_{k=0}^{n-1}|g(\lambda_{n,k})|^q .
 \]
 Now let $A_n = \{ z\in \DD \SUCHTHAT 1-\frac{3}{2n} \le |z| \le 1-\frac{1}{2n} \}$. Since
 $|g|^q$ is subharmonic,
 \[
   | g(\lambda_{n,k})|^q \le  \frac{4 n^2}{\pi} \int_{B(\lambda_{n,k}, \frac{1}{2n})}  | g(z) |^q dz
 \]
 And observing that the discs are disjoint subsets of $A_n$  yields
 \begin{align*}
   \frac1{n} \sum_{k=0}^{n-1}  |g(\lambda_{n,k})|^q
   \le  & \;  \frac{4 n}{\pi}  \int_{A_n} | g(z) |^q dz\\
   = & \; \frac{4 n}{\pi} \int_{ 1-\frac{3}{2n}}^ {1-\frac{1}{2n}} \int_0^{2\pi}|g(re^{it})|^q\,r\,dt\,dr
   \le 8\norm{g}_{H^q}^q. 
 \end{align*}
 This shows assumption \ref{item:Bonsall-Z-upper} of
 Theorem~\ref{thm:bonsall}.
 \proofskip
 For the converse estimate, we first consider the dense subspace
 $C(\overline{\DD})\cap H^q$.  For such functions, uniform continuity
 on $\overline{\DD}$ and the classical convergence of Riemann sums
 imply
\[
  \frac1n \sum_{k=0}^{n-1} |g(\lambda_{n,k})|^q \longrightarrow   \int_{\TT} |g(\zeta)|^q\,dm(\zeta).
\] Hence
\[
  \norm{g}_{H^q}\lesssim
\norm{ (\langle g,\varphi(\lambda)\rangle)_\lambda }_{Z'}.
\]
Since the sampling operator is bounded by the upper estimate already
established, the lower estimate extends to all of $H^q$ by density.
\end{proof}

Extending the preceding construction to more general radial
configurations leads to results in \cite{SperanskyTerekhin:2019} and
\cite{BaranovBatenev:2023}, which we recall without proof.

\begin{theorem}\label{thm:BB}[Speransky--Terekhin and Baranov--Batenev]
  Let $1<p<\infty$, let $r_n\uparrow 1$, and let $(m_n)$ be a sequence
  of positive integers. Consider the radial blocks
  \[
    \Lambda_n   =\bigl\{   r_n e^{2\pi i k/m_n}: k=0,\ldots,m_n-1      \bigr\},
    \qquad
    \Lambda=\bigsqcup_{n\ge1}\Lambda_n.
  \]
  If $\inf_{n\ge1}m_n(1-r_n)>0$, then the family
  $(k_\lambda)_{\lambda\in\Lambda}$ is a representing system for
  $H^p$.  If, moreover, $ m_n(1-r_n)\asymp 1$, then, with
  \[
    Z_\Lambda       :=\ell^1\bigl(\ell^p(\Lambda_n)\bigr),
    \qquad
    \varphi(\lambda)       :=\frac{k_\lambda}{\norm{ k_\lambda  }_{H^p}},
  \]
  every $f\in H^p$ admits a representation
  \[
    f=\sum_{\lambda\in\Lambda}        z_\lambda \, \varphi(\lambda),    \qquad z\in Z_\Lambda,
  \]
  satisfying
  \[
    \inf \Bigl\{ \norm{z}_{Z_\Lambda} \SUCHTHAT z\in Z_\Lambda \;
    \text{with}\; Sz=f \Bigr\} \; \asymp \; \norm{ f }_{H^p}.
  \]
\end{theorem}

The second assertion remains valid for more general configurations
$\lambda_{n,k}=r_ne^{i\theta_{n,k}}$, provided the successive angular
gaps are uniformly comparable to $1/m_n$, see
\cite[Section~3]{BaranovBatenev:2023}.  The proof in our framework,
even of a slight generalisation by the use of radial blockings,
i.e. sequences in annuli defined by increasing radii $r_n \uparrow 1$
would reproduce large parts of the arguments of the cited
authors, so we do not pursue it here.

\subsection{Necessity of quotient representations}\label{ss:quotient-example}
The example of section~\ref{ss:ST-example} can be re-used to
illustrate why the quotient space construction in
Theorem~\ref{thm:bonsall} is actually necessary.  Let
$r_n = 1- \frac1n$ and put
\[
  \psi_n=n^{-\sfrac1q} k_{r_n}.
\]
For each $n$, let $w_n\in Z$ be supported at the coordinate
corresponding to $r_n$, with value
$n^{-\sfrac1q}\norm{k_{r_n}}_{H^p}$. Then
\[
  S w_n=\psi_n,\qquad \norm{w_n}_Z\asymp1.
\]
Now let $a_N=\sum_{n=N}^{2N-1}(w_{2n}-w_{2n+1})$. Since all
coordinates live in distinct rings, $\norm{ a_N }_Z \asymp N$. But
close kernels nearly cancel. More precisely,
\[
  \psi_j-\psi_{j+1}
  =\bigl(j^{-\sfrac1q}-(j+1)^{-\sfrac1q}\bigr)k_{r_j}
   +(j+1)^{-\sfrac1q}(k_{r_j}-k_{r_{j+1}})
\]
By Lemma~\ref{lem:Hp-norms} the first term is $\calO( 1/j )$. The
second term is controlled by resolvent identity,
\[
  k_{r_j} - k_{r_{j+1}} = \frac{(r_j - r_{j+1})z }{ (1- r_{j}z)(1-r_{j+1}z) }
\]
so using Hölder's inequality and Lemma~\ref{lem:Hp-norms} with exponent $2p$, 
\[
  (j+1)^{-\sfrac1q}\norm{k_{r_j}-k_{r_{j+1}}}_{H^p}
  \lesssim j^{-\sfrac1q - 2} \norm{  k_{r_j} }_{H^{2p}} \norm{  k_{r_{j+1}} }_{H^{2p}}
  \lesssim 1/j.
\]
This means  that
\[
  \bignorm{ S a_N }_{H^p} \le \sum_{n=N}^{2N-1} \bignorm{ \psi_{2n} - \psi_{2n+1} }_{H^p}
  \lesssim  \sum_{n=N}^{2N-1} \frac1n = O(1)
\]
By the already established quotient theorem, there exists another
coefficient vector $c_N$ with $S c_N = S a_N$ and
$\norm{ c_N }_Z = O(1)$.  Thus $a_N-c_N \in \ker(S)$, and we have an
arbitrarily large discrepancy between two coefficient norms.

\subsection{Representations with zero density}

\begin{figure}[h]
  \centering
\begin{tikzpicture}[scale=8]
  \def\A{2}
  \pgfmathsetmacro{\tmax}{acos(1/\A)}
  \pgfmathdeclarefunction{rr}{1}{%
    \pgfmathparse{2*\A*(\A*cos(#1)-1)/(\A*\A-1)}%
  }
\begin{scope}
  \clip (-0.6,-0.1) rectangle (0.6,0.2);
  \fill[blue!15]
    plot[domain=-\tmax:\tmax,samples=150,variable=\t]
      ({-rr(\t)*sin(\t)},{rr(\t)*cos(\t)})
    -- cycle;
  \draw[dotted] (0,1) circle (1);
  \draw[blue,thick]
    plot[domain=-\tmax:\tmax,samples=150,variable=\t]
      ({-rr(\t)*sin(\t)},{rr(\t)*cos(\t)});
  \fill (0,0) circle (0.012);
  \node[below] at (0,0) {$\zeta$};
  \draw (0, 0.1) node {\small $\Gamma_A(\zeta)$};
  \def\lx{0.12}
  \def\ly{0.10}
  \pgfmathsetmacro{\rho}{sqrt(\lx*\lx+(\ly-1)^2)}
  \pgfmathsetmacro{\tc}{atan2(\lx,1-\ly)}
  \pgfmathsetmacro{\beta}{
    acos(1-(\A*\A-1)*(1-\rho)^2/(2*\rho))
  }
  \draw[very thick]
    plot[domain=\tc-\beta:\tc+\beta,
         samples=100,variable=\t]
      ({sin(\t)},{1-cos(\t)});
  \fill (\lx,\ly) circle (0.008);
  \draw (\lx+0.05, \ly) node {\small $\lambda$};    
  \draw[densely dotted]
    (\lx,\ly) -- ({sin(\tc)},{1-cos(\tc)});
  \draw ({sin(\tc)},{1-cos(\tc)-0.05}) node {\small $I_\lambda$};
\end{scope}
\end{tikzpicture}
\caption{The shadow $I_\lambda$ of a point and the outer part of the Stolz
  region $\Gamma_A(\zeta)$ on a boundary vertex.}
  \label{fig:shadow-stolz}
\end{figure}
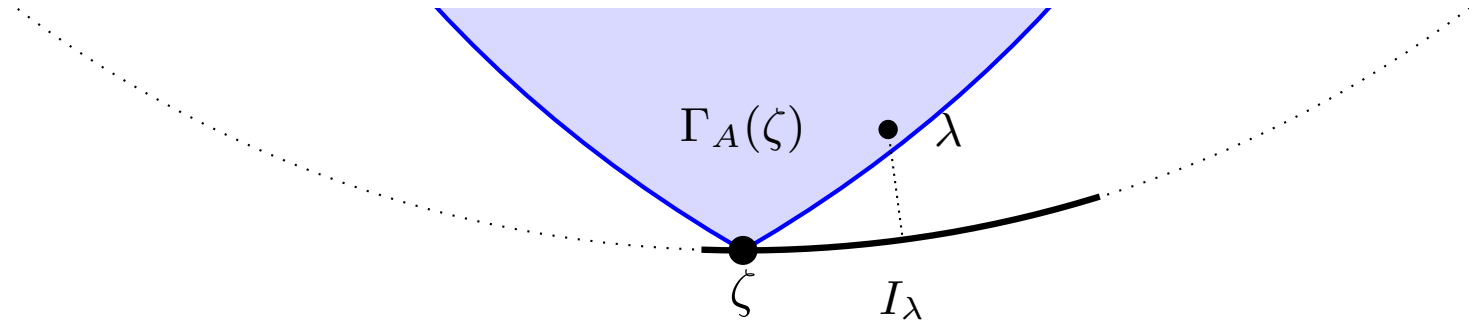

We now come to the first set of main results. In
\cite{BaranovBatenev:2023}, the authors ask whether there exist
representing systems of Cauchy kernels in $H^p$, $1<p<\infty$, with
vanishing upper density
\[
  D^+(\Lambda)   := \limsup_{r\uparrow1}   \quad (1-r) \,\CM( \bigl\{\lambda\in\Lambda \SUCHTHAT |\lambda|<r  \bigr\} )
\]
where $\kappa$ means the counting measure. We shall give an
affirmative answer. The idea is to abandon radial blockings altogether
and to use non-tangential approach regions instead. Fix an aperture
parameter $A>1$. For $\zeta\in \TT$, define the Stolz domain with vertex
$\zeta$ and aperture $A$ as
\[
  \Gamma_A(\zeta)   :=\{z\in\DD \SUCHTHAT |z-\zeta| < A (1-|z|) \}.
\]
Conversely, each $\lambda\in\DD$ determines its \emph{shadow}
\[
  I_\lambda   :=\{\zeta\in\TT \SUCHTHAT \lambda\in\Gamma_A(\zeta)\},
\]
see Figure~\ref{fig:shadow-stolz}.

\noindent For a countable set $\Lambda\subset\DD$, put
$\Lambda(\zeta):=\Lambda\cap\Gamma_A(\zeta)$ and let
\[
  \MF_\Lambda g(\zeta) :=\sup_{\lambda\in\Lambda(\zeta)}|g(\lambda)|
\]
be the sampled maximal function (where the supremum over the empty set
is understood to be zero).  Clearly, $\MF_\Lambda g\le \MF_\Gamma g$,
where $\MF_\Gamma$ denotes the classical non-tangential maximal
function.  Consequently,
\begin{equation}  \label{eq:sampled-max-upper}
  \norm{ \MF_\Lambda g  }_{L^q(\TT)}   \le C_{A,q}\norm{ g  }_{H^q},   \qquad 1<q<\infty.
\end{equation}
Such maximal functions were considered, for example in
\cite{BrunaNicolauOyma}; the corresponding sampling theory was studied
by Thomas \cite{Thomas:Sampling1998}, see also \cite{MassanedaThomas}.
We recall the relevant geometric criterion.

\begin{definition}
  We say that $\Lambda$ is \emph{non-tangentially dense almost   everywhere} if, for some $A>1$,
  \begin{equation}    \label{eq:nt-density}
    \lim_{r\to 1-}  \;    m\left(  \bigcup_{\lambda\in\Lambda, |\lambda|>r}   I_\lambda  \right)=1.
  \end{equation}
\end{definition}

\noindent Since the sets $U_r(\Lambda) :=\bigcup_{\lambda\in\Lambda, |\lambda|>r}  I_\lambda$ 
decrease as $r\uparrow1$, continuity of measure from the above shows that
condition \eqref{eq:nt-density} is equivalent to 
\[
   E_\Lambda  :=\bigcap_{j\ge2} U_{1-1/j}(\Lambda)
 \]
 having full measure $m(E_\Lambda)=1$.  In particular, almost every
 boundary point is approached non-tangentially (with uniform aperture)
 by points of $\Lambda$.  The following theorem is due to Thomas
 \cite[Theorem~1]{Thomas:Sampling1998} building upon an earlier result in 
 Bruna-Nicolau-\O{}yma \cite{BrunaNicolauOyma}.

\begin{theorem}[Thomas] \label{thm:thomas:sampling}
  Let $\Lambda\subset\DD$ be countable and let $A>1$.  Then the
  following assertions are equivalent:
  \begin{enumerate}
  \item Almost every $\zeta\in\TT$ is a non-tangential accumulation
    point of $\Lambda$ within $\Gamma_A(\zeta)$, that is,
      \[
        m(E_\Lambda)=1.
      \]
    \item For some (equivalently, every) $0<q<\infty$ there exists
      $c_q>0$ such that
      \[
        \norm{ \MF_\Lambda g }_{L^q(\TT)} \ge c_q \norm{ g }_{H^q}, \qquad g\in H^q.
      \]
  \end{enumerate}
  In this case, one actually has
  \begin{equation} \label{eq:sampled-max-equivalence}
    \norm{ g  }_{H^q}     \le \norm{ \MF_\Lambda g  }_{L^q(\TT)}     \le C_{A,q}\norm{ g  }_{H^q},     \qquad g\in H^q.
  \end{equation}
\end{theorem}

We now combine Thomas's sampling theorem with our abstract quotient
principle to obtain a characterization of kernel representations in
terms of the geometry of their shadows.  Observe that, in contrast to
Theorem~\ref{thm:BB}, we work here with non-normalized kernels.

\begin{theorem}\label{thm:nt-representation}
  Let $1<p<\infty$, let $\Lambda\subset\DD$ be countable, and fix an
  aperture $A>1$.  For a finitely supported coefficient sequence
  $(c_\lambda)$, let
  \[
    F_c(\zeta)     := \sum_{\lambda\in\Lambda}       |c_\lambda| \, \frac{\eins_{I_\lambda}(\zeta)}{m(I_\lambda)}.
  \]
  Define
  \[
    \CS_\Lambda     :=\overline{c_{00}(\Lambda)}^{\norm{\cdot}_{\CS_\Lambda}},
    \quad \text{where}\quad   \norm{c}_{\CS_\Lambda} := \norm{F_c}_{L^p(\TT)}.
  \]
  Then the synthesis operator
  \[
    S:\CS_\Lambda\longrightarrow H^p,     \qquad     Sc=\sum_{\lambda\in\Lambda}c_\lambda k_\lambda,
  \]
  is bounded, with unconditional convergence of the series.  Moreover,
  $S$ is surjective if and only if $\Lambda$ is non-tangentially dense
  almost everywhere with respect to $\Gamma_A$.  In this case,
  $(k_\lambda)_{\lambda\in\Lambda}$ is an unconditionally representing
  system for $H^p$, and
  \[
    H^p\simeq\CS_\Lambda/\ker S.
  \]
\end{theorem}

\begin{remark}
  Since removing finitely many points preserves non-tangential
  density, Theorem~\ref{thm:nt-representation} shows that every
  individual kernel can be represented using the remaining ones. In
  particular, the synthesis map has a nontrivial kernel, giving
  another concrete motivation for the quotient construction.
\end{remark}

\begin{proof}
  Let $X=H^p$, and denote by $q$ the conjugate exponent.  Recall that
  $X'=(H^p)'$ is identified with $H^q$, with equivalence of norms,
  through the pairing
  \[
    \langle \ell_g,f\rangle
    =\int_{\TT} f(\zeta)\overline{g(\zeta)}\,dm(\zeta), \qquad g\in H^q,\quad f\in     H^p.
  \]
  In particular,
  $\langle\ell_g,k_\lambda\rangle =\overline{g(\lambda)}$.  We first
  show that the synthesis operator $S$ is bounded, without any
  assumption on the distribution of $\Lambda$.  For
  $c\in c_{00}(\Lambda)$ and $g\in H^q$, we have
\begin{align*}
  |\langle\ell_g, Sc \rangle|
  & \le  \; \sum_{\lambda\in\Lambda}|c_\lambda|\,|g(\lambda)|
    \; = \; \int_\TT     \sum_{\lambda\in\Lambda}    |c_\lambda| \frac{\eins_{I_\lambda}(\zeta) }{m(I_\lambda)}    |g(\lambda)|    \,dm(\zeta)\\
  & \le  \; \int_\TT F_c(\zeta) (\MF_\Gamma g)(\zeta)\,dm(\zeta)
    \; \le \; \norm{c}_{\CS_\Lambda}     \norm{\MF_\Gamma g}_{L^q(\TT)} \\
  & \lesssim_{A,q}  \; \norm{c}_{\CS_\Lambda}    \norm{g}_{H^q}.
\end{align*}
By duality, $S$ extends to a bounded operator
$S:\CS_\Lambda\to H^p$.  Moreover, coordinate projections are
contractive on $\CS_\Lambda$ and finite sequences are dense.
Consequently, the series defining $Sc$ converges unconditionally in
$H^p$ for every $c\in\CS_\Lambda$.

\medskip \noindent It remains to characterize the surjectivity of $S$.
First assume that $\Lambda$ is non-tangentially dense almost
everywhere. We then consider the measure space
  \[
    \Omega :=\{(\zeta,\lambda)\in\TT\times\Lambda \SUCHTHAT \zeta\in  I_\lambda\},
  \]
  equipped with the restriction of $m\otimes\CM$, where $\CM$ denotes
  counting measure.  Define the Banach function space
  \[
    Z:=L^p(\TT;\ell^1(\Lambda(\zeta))),
  \]
  more explicitly, by
  \begin{equation}    \label{eq:nt-coefficient-norm}
    \norm{ a  }_Z     :=    \Bignorm{    \sum_{\lambda\in\Lambda(\cdot)}        |a(\cdot,\lambda)|       }_{L^p(\TT)}.
  \end{equation}
  This space has absolutely continuous norm.  Its Köthe dual is
  given by
  \[
    \norm{ b }_{Z'} = \Bignorm{ \sup_{\lambda\in\Lambda(\cdot)}       |b(\cdot,\lambda)|  }_{L^q(\TT)}.
  \]
  This follows from the duality between $\ell^1$ and $\ell^\infty$,
  followed by $L^p$--$L^q$ duality.  Set $\varphi(\zeta,\lambda):=k_\lambda\in H^p$.
  Since $\Lambda$ is countable, $\varphi$ is strongly measurable.
  The corresponding analysis operator   $T:X'\to Z'$ is
  \[
    (T\ell_g)(\zeta,\lambda)
      =\langle\ell_g,\varphi(\zeta,\lambda)\rangle
      =\overline{g(\lambda)}.
  \]
  Consequently, using Theorem~\ref{thm:thomas:sampling}, 
  \[
    \norm{ T\ell_g }_{Z'} \; = \; \norm{ \MF_\Lambda g }_{L^q(\TT)} \; \overset{  \eqref{eq:sampled-max-equivalence}}{\asymp} \; \norm{ g }_{H^q}  \;  \asymp \;   \norm{ \ell_g }_{X'}.
  \]
  Thus $T$ is an isomorphic embedding.  By Theorem~\ref{thm:bonsall},
  its adjoint restricts to a bounded surjection
  \[
    S_0:=T'|_Z:Z\longrightarrow H^p.
  \]
  More explicitly, for $a\in Z$ and $g\in H^q$,
  \[
    \langle\ell_g, S_0 a\rangle
    =\langle T\ell_g,a\rangle =\sum_{\lambda\in\Lambda}      \int_{I_\lambda}      a(\zeta,\lambda)\overline{g(\lambda)}      \,dm(\zeta).
  \]
  Setting 
  \[
    c_\lambda := \int_{I_\lambda}a(\zeta,\lambda)\,dm(\zeta),
  \]
  for $\lambda \in \Lambda$ we therefore obtain
  $S_0a=\sum_{\lambda\in\Lambda} c_\lambda k_\lambda$.   Let
  $\pi:Z\to\mathcal Q$ be the coefficient map $a\mapsto(c_\lambda)$,
  where $\mathcal Q=\pi(Z)$ is equipped with the quotient norm. Since
  $\ker\pi\subseteq\ker S_0$, the operator $S_0$ factors through
  $\mathcal Q$, giving a bounded surjection $S:\mathcal Q\to H^p$ with
  $S_0=S\circ\pi$.
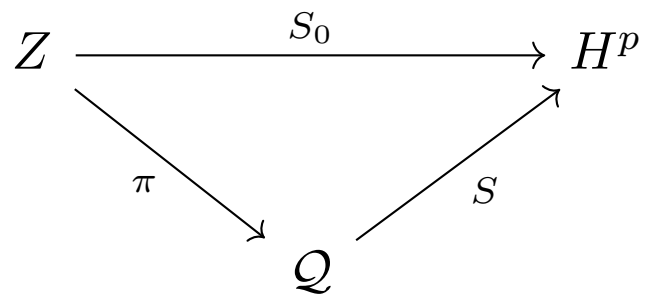
\[
\begin{tikzcd}[column sep=large, row sep=large]
Z \arrow[rr,"S_0"] \arrow[rd,,"\pi"']
  & & H^p \\
& \mathcal Q \arrow[ur,"S"'] &
\end{tikzcd}
\]
We claim that $\mathcal Q$ identifies isomorphically with the space
$\CS_\Lambda$ defined in the statement. Notice that this
identification does not require non-tangential density.  Indeed, for a
given sequence $(c_\lambda)$, let
\[
  a(\zeta, \lambda) := c_\lambda \frac{ \eins_{I_\lambda} } {m(I_\lambda)}
\]
Then $\pi(a)= c$ and $\norm{a}_Z = \norm{ F_c }_{L^p}$. Moreover, for any $a \in Z$ such that $\pi(a) = c$,
and positive function $h \in L^q(\TT)$,
\[
  \int_\TT F_c \, h \,dm = \sum_{\lambda \in \Lambda} |c_\lambda| \langle h \rangle_{I_\lambda}
\]
where $\langle h \rangle_I$ is the average of $h$ on $I$. It follows that, denoting 
the uncentered Hardy--Littlewood maximal operator as $\MF$, 
\begin{align*}
  \int_\TT F_c \, h \,dm
  =   & \; \sum_{\lambda \in \Lambda} \Bigl|  \int_{I_\lambda} a(\zeta, \lambda)\,dm(\zeta) \Bigr| \langle h \rangle_{I_\lambda}\\
  \le & \; \sum_{\lambda \in \Lambda}   \int_{I_\lambda} \bigl| a(\zeta, \lambda)\bigr|  \,dm(\zeta)  \langle h \rangle_{I_\lambda}\\
  \le & \; \int_\TT  \sum_{\lambda \in \Lambda(\zeta)}   \bigl| a(\zeta, \lambda)\bigr|   (\MF h)(\zeta)\,dm(\zeta)  \\
  \le & \; C_q \norm{a}_Z \norm{h}_{L^q(\TT)}  
\end{align*}
Taking supremum over all such $h$ in the unit ball yields
$\norm{ F_c }_{L^p} \le C_q \norm{ a }_Z$, and taking the infimum over
all $a\in Z$ satisfying $\pi(a) = c$ finally achieves
\[
  \bignorm{ F_c}_{L^p(\TT)} \le C_q \bignorm{(c_\lambda)}_{\mathcal Q} \le C_q \bignorm{ F_c }_{L^p(\TT)}.
\]
Since finite-coordinate truncations converge in $Z$, the finitely
supported sequences are dense in $\mathcal Q$. The preceding
inequalities show that the two norms are equivalent on
$c_{00}(\Lambda)$.  Consequently, their completions are isomorphic,
and hence
\[
  \mathcal Q\simeq\CS_\Lambda.
\]
Under this identification, $S$ is precisely the synthesis operator
appearing in the statement.  The quotient principle now gives
\[
  H^p\simeq\CS_\Lambda/\ker S.
\]
For the converse implication, suppose that $S:\CS_\Lambda\to H^p$ is
surjective. Then its adjoint $S'$ is bounded below. Recall that the
coefficient-space identification of $\CS_\Lambda$ with $\mathcal Q$
applies without any non-tangential density hypothesis. But the adjoint
norm is equivalent to the norm of the analysis operator $T$, namely
\[
  \norm{T\ell_g}_{Z'}   =\norm{\MF_\Lambda g}_{L^q(\TT)}.
\]
Consequently,
\[
  \norm{g}_{H^q}   \lesssim \norm{\MF_\Lambda g}_{L^q(\TT)}.
\]
By Theorem~\ref{thm:thomas:sampling}, $\Lambda$ is non-tangentially dense
almost everywhere.
\end{proof}

\begin{remark}
  Theorem~\ref{thm:nt-representation} contains the representing-system
  assertion of Theorem~\ref{thm:BB} as a special case.  Indeed,
  suppose that
  \[
    \delta := \inf_n m_n(1-r_n)>0.
  \]
  For every $\zeta\in\TT$ and every $n$, there exists
  $\lambda\in\Lambda_n$ such that
  \[
    |\zeta-\lambda| \le (1-r_n)+\frac{\pi}{m_n} \le    \left(1+\frac{\pi}{\delta}\right)(1-r_n).
  \]
  Hence, choosing $A>1+\pi/\delta$, we have
  \[
    \bigcup_{\lambda\in\Lambda_n}I_\lambda=\TT \qquad\text{for every
    }n,
  \]
  and condition~\eqref{eq:nt-density} follows.  Notice, however, that
  the coefficient space $\CS_\Lambda$ constructed here differs
  from the block space $\ell^1(\ell^p(\Lambda_n))$ appearing in
  Theorem~\ref{thm:BB}, even in the critical case
  $m_n(1-r_n) \asymp 1$. In this case, we have shadows with uniformly bounded overlap for every fixed aperture. Writing
  \[
    F_n(\zeta)  = \sum_{\lambda \in \Lambda_n} \frac{ |z_\lambda|  \eins_{I_\lambda}(\zeta) }{\norm{k_\lambda}_{H^p}\,m(I_\lambda) },
  \]
  we have therefore
  \[
    \bignorm{ F_n }_{L^p}^p \asymp \sum_{\lambda \in \Lambda_n}  |z_\lambda|^p
  \]
  so that Theorem~\ref{thm:BB} measures an $\ell_1( L^p)$-norm,
  whereas Theorem~\ref{thm:nt-representation} measures a
  $L^p( \ell_1)$-norm.  This looks surprising at first glance, but can
  be convincingly explained in the dual spaces: where
  Theorem~\ref{thm:BB} controls
  \[
    \sup_n \Bigl( \sum_{\lambda \in \Lambda_n} |g(\lambda)|^q \Bigr)^{\sfrac1q},
  \]
  Theorem~\ref{thm:nt-representation} rather controls
  \[
    \Bignorm{  \zeta \mapsto \sup_{\lambda \in \Lambda\cap\Gamma_A(\zeta)} |g(\lambda)|  }_{L^q}.
  \]
  These are very different sampling norms. Nevertheless, both are
  equivalent to $\norm{ g}_{H^q}$ --- this equivalence being asserted
  only for sequences of values of \emph{analytic functions}, not for
  arbitrary scalar sequences on $\Lambda$.
\end{remark}

Our new coefficient geometry, as we will see now, allows to answe the
density question raised in \cite{BaranovBatenev:2023}. The
construction is inspired by the random sampling argument of Massaneda
and Thomas \cite[Theorem 3]{MassanedaThomas}. Here, however, the
points on each circle are equally spaced, and only the circles
themselves are independently rotated.

\begin{coro}\label{cor:zero-density}
  There exists a sequence $\Lambda\subset\DD$ satisfying
  $D^+(\Lambda)=0$ such that $(k_\lambda)_{\lambda\in\Lambda}$ is an
  unconditionally representing system for $H^p$ for every
  $1<p<\infty$.
\end{coro}
\begin{proof}
  Let $r_n=1-2^{-n}$ and $m_n=\left\lfloor\frac{2^n}{n}\right\rfloor$
  where $n\ge2$.  Let $(\theta_n)$ be independent random variables,
  uniformly distributed on $[0,2\pi)$, and put
  \[
    \Lambda_n     :=    \left\{      r_n\exp\left(        i\theta_n+\frac{2\pi ik}{m_n}
      \right) \SUCHTHAT       k=0,\ldots,m_n-1     \right\},
    \qquad
    \Lambda=\bigcup_{n\ge2}\Lambda_n.
  \]
  First, if $r_n\le r<r_{n+1}$, then
  \[
    (1-r)\; \CM \{\lambda\in\Lambda \SUCHTHAT|\lambda|<r\}
    \le 2^{-n}\sum_{k=2}^{n}m_k     \lesssim\frac1n.
  \]
  Hence $D^+(\Lambda)=0$, independently of the random rotations.  It
  remains to verify the non-tangential density condition. Fix an
  aperture $A>1$ and write
  \[
    V_n:=\bigcup_{\lambda\in\Lambda_n}I_\lambda.
  \]
  For $\lambda\in\Lambda_n$, the normalized length of its shadow
  satisfies
  \[
    m(I_\lambda)\asymp_A (1-r_n) =  2^{-n}.
  \]
  Moreover, since $m_n2^{-n}\to0$, the shadows belonging to a single
  ring $\Lambda_n$ become pairwise disjoint for all sufficiently large
  $n$.  Therefore
  \[
    m(V_n)\asymp_A m_n2^{-n} \asymp\frac1n.
  \]
  Fix $\zeta\in\TT$. By rotational invariance,
  \[
    \mathbb P(\zeta\in V_n)=m(V_n)\asymp\frac1n.
  \]
  The events $\{\zeta\in V_n\}$ are independent because the rotations
  $\theta_n$ are independent. The second Borel--Cantelli lemma
  therefore yields
  \[
    \PP(\zeta\in V_n \text{ for infinitely many }n)=1, 
  \]
  i.e. $\zeta \in \limsup_{n} V_n$ almost surely.  Integrating over
  $\zeta\in\TT$ and using Tonelli's theorem, we conclude that
\[
      \EE\left[m\left(\limsup_{n\to\infty} V_n\right)\right]
  =   \int_{\TT}     \PP\left(\zeta\in\limsup_{n\to\infty}V_n\right)    \,dm(\zeta)
  =   1.
\]
Since $0\le m(\limsup_n V_n)\le1$, we conclude that
\[
  m\left(\limsup_{n\to\infty} V_n \right)=1   \qquad\text{almost surely}.
\]
Since $r_n\uparrow1$, this is precisely condition
\eqref{eq:nt-density}.  Choosing any realization for which this holds,
Theorem~\ref{thm:nt-representation} gives the assertion.
\end{proof}

\section{Random kernel points}

We now apply the shadow criterion to two models of random points in
$\DD$.

\subsection{Independent random angles}

Let $(\theta_n)$ be independent uniform random variables on
$[0, 2\pi)$, $(r_n)$ a nondecreasing deterministic sequence with $r_n\to1$, and
$\lambda_n = r_n e^{i \theta_n}$. We write
$\Lambda(\om) = \{ \lambda_n(\om)\SUCHTHAT n \ge 1 \}$. For a fixed
$\zeta \in \TT$, write
\[
  E_n(\zeta) = \{ \om \in \Om\SUCHTHAT \zeta \in I_{\lambda_n(\om)} \}
\]
Since the angles are uniform,
\[
  \PP( E_n(\zeta))  = m(I_{r_n})  \; \asymp_A \; (1-r_n).
\]
The events $E_n(\zeta)$ are independent, and so Borel-Cantelli gives
\[
  \PP( \zeta \in I_{\lambda_n} \;\text{infinitely often}\; )
  =
  \left\{
    \begin{array}{lcl}
      1 & \text{if} & \sum_n  (1-r_n) = +\infty \\
      0 & \text{if} & \sum_n  (1-r_n) < +\infty 
    \end{array}
  \right.
\]
In the divergent case, integrating in $\zeta$ and Tonelli's theorem imply
\[
  \EE \bigl[   m( \limsup I_{\lambda_n} ) \bigr]
  = \int_\TT \PP(  \zeta \in  \limsup I_{\lambda_n} )\,dm(\zeta) = 1
\]
and since $m$ takes values in $[0,1]$ this implies
$m( \limsup I_{\lambda_n})=1$ almost surely, i.e. non-tangential
density. By Theorem~\ref{thm:nt-representation},
$\{ k_{\lambda_n} \SUCHTHAT n\ge 1 \}$ is almost surely representing.
In the convergent case, $ m( \limsup I_{\lambda_n} )=0$ almost surely
for the same reason. This is precisely the random sampling criterion
proved by Massaneda--Thomas \cite[Theorem 3]{MassanedaThomas} in the
unit-disc case.

\begin{coro}\label{cor:random-angles}
  Let $(r_n)$ be a nondecreasing deterministic sequence of radii
  converging to $1$, and let $\lambda_n$ be obtained by independent,
  uniformly distributed random angles. Then the reproducing kernels
  form a representing system almost surely if and only if
  \[
    \sum_n (1-r_n) = +\infty.
  \]
\end{coro}

 \subsection{Poisson point processes}
 Finally let us consider a locally finite, rotationally invariant
 Poisson point process with intensity $d\mu(re^{it}) = d\nu(r) dm(t)$.
 Thus $\Lambda$ is a random point set characterized by
 \begin{enumerate}[label=(\alph*), ref=(\alph*),leftmargin=*]
 \item{}(Poisson counts). For measurable $B$ with $\mu(B)<\infty$,
   \[
     \PP( \#(\Lambda \cap B)=k ) = e^{-\mu(B)} \frac{\mu(B)^k}{k!}
   \]
 \item{}(Independence). For pairwise disjoint measurable sets
   $B_1,\ldots,B_n$, the random counts
   $N(B_i)=\#(B_i\cap\Lambda)$ are independent.
 \item{}(Local finiteness). The radial intensity $\nu$ satisfies
   $\nu([0,R])<\infty$ for every $R<1$.
 \end{enumerate}
 The intensity is rotationally invariant. By Campbell's
 formula, for positive measurable $h: \DD \to [0, \infty)$,
 \[
   \EE\Bigl[ \sum_{\lambda \in \Lambda} h(\lambda) \Bigr] = \int_\DD h(z) \,d\mu(z)
 \]
 Choosing $h(z) = 1-|z|$, an immediate consequence is that if the
 integral
 \begin{equation}   \label{eq:integral}
        \int_0^1 (1-r)\,d\nu(r)
 \end{equation}
 is finite, then
 \[
      \sum_{\lambda \in \Lambda} (1-|\lambda|) < \infty
 \] 
 almost surely. Thus the random points form a Blaschke sequence almost
 surely, and cannot give a representing system for $H^p$.  Now
 consider the case of divergence of \eqref{eq:integral}. Recall the
 notation
 \[
   E_\Lambda = \bigcap_{j\ge 2} \bigcup_{\lambda \in \Lambda, |\lambda| > 1-\sfrac1j}  I_\lambda,
 \]
 i.e. the points of $\TT$ that can be approached non-tangentially (with
 uniform aperture $A$) from within $\Lambda$.  Looking at any outer
 part of a Stolz region
 \[
   B_R(\zeta) = \Gamma_A(\zeta) \cap \{ |z|>R \},
 \]
 we calculate its intensity
 \[
   \mu(B_R) = \int_R^1 m(I_r)\,d\nu(r) \;\asymp\;  \int_R^1 (1-r) \,d\nu(r) = +\infty
 \]
 so that, for each fixed $\zeta$, every $B_{1-1/j}(\zeta)$
contains infinitely many Poisson points almost surely. Indeed, this
follows by exhaustion with sets of finite intensity. Taking the
countable intersection over $j\ge2$ gives
$\PP(\zeta\in E_\Lambda)=1$ for every fixed $\zeta$.
 Then the Tonelli argument gives
 \[
   \EE \bigl[ m(E_\Lambda) \bigr] = \int_\TT \PP( \zeta \in E_\Lambda )\,dm(\zeta) = 1
 \]
 and we infer $m(E_\Lambda)=1$ almost surely. Now by Theorems
 \ref{thm:thomas:sampling} and \ref{thm:nt-representation} we obtain
 the following result.

 \begin{theorem}\label{thm:poisson}
   Let $1<p<\infty$ and consider a locally finite, rotationally
   invariant Poisson point process $\Lambda$ with intensity
   $d\mu(re^{it}) = d\nu(r) dm(t)$ on $\DD$. Then the set
   \[
     \{ k_\lambda \SUCHTHAT \lambda \in \Lambda \}
   \]
   is almost surely representing in $H^p$ if and only if the integral
   \eqref{eq:integral} diverges.  In that case, the representation is
   realized through the shadow coefficient space of
   Theorem~\ref{thm:nt-representation}.
 \end{theorem}

 \section{Explicit representations}\label{sec:explicit}

 While Theorem~\ref{thm:nt-representation} characterizes
 representations through the shadow coefficient space, it does not
 provide explicit coefficients. We now turn to approximations
 constructed directly from the sampled values of a function.  Related
 reconstruction problems have been studied by Cima and Stessin
 \cite{CimaStressin}, who use Gram matrices and orthogonal
 projections.  Sampling inequalities for polynomial subspaces of Hardy
 spaces have also been investigated by Gröchenig and Ortega-Cerdà
 \cite{Grochenig-OrtegaCerda}. Our approach instead relies on
 geometric quadrature formulas and non-tangential maximal estimates.
 Recall that non-tangential density assumption says that, for every
 $r<1$,
 \[
   m\Bigl( \bigcup_{\lambda \in \Lambda, |\lambda|>r}  I_\lambda \Bigr) = 1
 \]
 For each $n\ge2$, we may therefore select finitely many points
$\lambda_{n,j}$, $j=1,\ldots,J_n$, with
$|\lambda_{n,j}|>1-\frac1n$, whose shadows cover a measurable set
$W_n$ of measure at least $1-2^{-n}$. Choose disjoint measurable parts
 $E_{n, j} \subseteq I_{\lambda_{n, j}} $ so that
 \[
   \bigcup_{j=1}^{J_n} E_{n,j}  =    \bigcup_{j=1}^{J_n} I_{\lambda_{n, j}} = W_n.
  \]
  For a given $H^p$ function $f$, put
  \begin{equation}   \label{eq:Qn}
    \begin{aligned}
      Q_n(f)(z) & := \sum_{j=1}^{J_n} R_n(f)_j  \; k_{\lambda_{n,j}}(z) \\
      R_n f & := c^{(n)}  : = \bigl( m(E_{n,j})\,f(\lambda_{n,j})   \bigr)_{j=1, \ldots, J_n}
    \end{aligned}
 \end{equation}
 We observe that the coefficient norm stays uniformly bounded. Indeed,
 let $h\ge 0$ be in $L^q$. Then
 \begin{align*}
         \int_\TT h(\zeta)F_{c^{(n)}}(\zeta)\,dm(\zeta)   
    =  & \;\;   \sum_{j=1}^{J_n}  \Bigl(\frac{1}{m(I_{\lambda_{n,j}})} \int_{I_{\lambda_{n,j}}}  h(\zeta) \,d\zeta \Bigr)  m(E_{n,j}) \, |f(\lambda_{n,j})| \\
   \le & \;\; \sum_{j=1}^{J_n}  \int_{E_{n,j}}  \MF_{\Gamma}(f) (\zeta)  \MF(h)(\zeta)\,dm(\zeta)\\
   \le & \;\; C_{p,A} \norm{f}_{H^p} \norm{h}_{L^q},
 \end{align*}
 and hence $\norm{ c^{(n)}(f) }_{\CS_\Lambda} \le C_{p,A} \norm{f}_{H^p}$.  Moreover, for $g \in H^q$,
 \begin{equation}   \label{eq:prepa-dom-cv}
   \dprod{ g }{Q_n f}
   =  \sum_{j=1}^{J_n} m(E_{n,j}) \, f(\lambda_{n,j})   \overline{g(\lambda_{n,j})}
   = \int_{W_n} f( \lambda_n(\zeta))  \overline{g(\lambda_{n}(\zeta)) }\,dm(\zeta)
 \end{equation}
 where $\lambda_n(\zeta)=\lambda_{n,j}$ for $\zeta\in E_{n,j}$.
This map is piecewise constant, hence measurable.  Since
 $E_{n,j} \subseteq I_{\lambda_{n,j}}$, every point
 $\zeta \in E_{n,j} $ satisfies $\lambda_{n,j} \in
 \Gamma(\zeta)$. It follows that $\lambda_n(\zeta) \to \zeta$ non-tangentially
 whenever $\zeta \in W_n$ for all large $n$. But by the first
 Borel--Cantelli lemma, applied to the deterministic sets
 $\TT \setminus W_n$ on the probability space $(\TT, m)$, almost every
 $\zeta$ belongs to $W_n$ eventually.  This implies that
 $f( \lambda_n(\zeta)) \to f^*(\zeta)$ almost everywhere, and
 likewise for $g$. Moreover,
 $|f(\lambda_{n,j})|\le \MF_\Gamma(f)(\zeta)$ for all
 $\zeta \in E_{n,j}$, so that the integrand in \eqref{eq:prepa-dom-cv}
 is pointwise dominated by the integrable function
 $\MF_\Gamma(f) \MF_\Gamma(g)$. By the dominated convergence theorem,
 \[
   \dprod{ g }{Q_n f} \quad \overset{n \to \infty}{\relbar\joinrel\relbar\joinrel\relbar\joinrel\to} \quad  \int_\TT f^*(\zeta)\overline{g^*(\zeta)}\,dm(\zeta).
 \]
 We have reached the following proposition.

 \begin{proposition}\label{prop:weak-convergence}
   Let $1<p<\infty$, $f \in H^p$ and $Q_n f$ be defined by
   \eqref{eq:Qn}. Then $Q_nf\rightharpoonup f$ weakly in $H^p$,
   \[
     \sup_n\norm{Q_n}_{\BOUNDED(H^p)}<\infty,
     \qquad
     \sup_n\norm{c^{(n)}(f)}_{\CS_\Lambda}
       \lesssim_{p,A}\norm{f}_{H^p}.
   \]
   More precisely, the operator bounds follow from the
   non-tangential maximal estimates in $H^p$ and $H^q$.
 \end{proposition}

 \begin{example}
   Let $f=\eins$, $r_n=1-2^{-n}$, and place $N=2^n$ equally
   spaced points on the circle of radius $r_n$. With
   $\omega=e^{2\pi i/N}$, the roots-of-unity identity in
   Lemma~\ref{lem:roots-of-unity} gives
   \begin{equation}\label{eq:Qn-von-1}
     Q_n \eins(z)=\frac1N\sum_{j=0}^{N-1}        \frac{1}{1-r_n\omega^{-j}z}       =\frac1{1-r_n^Nz^N}.
   \end{equation}
   Consequently,
   \[
     \norm{Q_n \eins -\eins }_{H^2}        =\frac{r_n^N}{\sqrt{1-r_n^{2N}}}       \longrightarrow\frac1{\sqrt{e^2-1}}>0.
   \]
   Thus weak convergence of the quadrature operators need not be
   strong convergence.
 \end{example}

 Since $(Q_n f)$ converges weakly, suitable convex combinations
 converge in norm by Mazur's theorem. In the setup of equally spaced
 points considered in the example above, a pleasant surprise happens.

 \begin{theorem}[Special Fej\'er phenomenon]
   \label{thm:fejer-std-config}
   Let $\Lambda$ consist of $n$ equally spaced points on each
   circle of radius $r_n=1-1/n$, $n\ge2$. For the quadrature
   operators $Q_n$ obtained from the equal angular sectors, put
   \[
     T_N=\frac1N\sum_{n=2}^{N+1}Q_n.
   \]
   Then $T_Nf\to f$ in $H^p$ for every $f\in H^p$ and $1<p<\infty$. The conclusion remains valid if each ring is
   rotated by an arbitrary angle $\theta_n$.
 \end{theorem}
 \begin{proof}
   We first take $\theta_n=0$. The equal sectors lie in the
   corresponding shadows once the aperture satisfies $A>1+\pi$. Hence
   Proposition~\ref{prop:weak-convergence} gives uniform boundedness
   of $Q_n$ on every $H^p$, $1<p<\infty$.

   \proofskip
   It suffices to prove strong convergence on polynomials. Fix
   $f(z)=z^m$. For $n>m$, Lemma~\ref{lem:roots-of-unity} yields
   \begin{equation}\label{eq:Q_n-von-z-hoch-m}
     Q_nz^m(z)=\frac{r_n^{2m}z^m}{1-r_n^nz^n},
   \end{equation}
   and hence
   \[
     Q_nz^m-z^m        =(r_n^{2m}-1)z^m         +\sum_{j\ge1}r_n^{nj+2m}z^{nj+m}.
   \]
   Since $1-r_n^{2m}\le 2m/n$, the Ces\`aro mean of the first
   term is $O_m(\log(N+2)/N)$. For each fixed $j$, the exponents
   $m+nj$ are distinct as $n$ varies. Orthogonality in $H^2$
   therefore implies
   \begin{align*}
     &\frac1N\Bignorm{\sum_{j\ge1} \;       \sum_{\max(1, m) < n\le N+1}     r_n^{nj+2m}z^{nj+m}}_{H^2}\\
     &\quad\le\frac1N\sum_{j\ge1} \;        \Bigl( \sum _{\max(1, m) < n\le N+1}               r_n^{2nj+4m}\Bigr)^{1/2} \\
     &\quad\le\frac1{\sqrt N}\sum_{j\ge1}e^{-j}    =\frac1{(e-1)\sqrt N}.
   \end{align*}
   The finitely many terms with $n\le m$ contribute
   $O_m(N^{-1})$ by uniform boundedness. Thus
   \[
     \norm{T_Nz^m-z^m}_{H^2}
       \lesssim_m\frac{\log(N+2)}N+\frac1{\sqrt N}\longrightarrow0.
   \]

   For a fixed polynomial $P$, the same lemma gives a uniform
   $H^\infty$ bound on $Q_nP$, since $r_n^n\le e^{-1}$.
   Hence $T_NP-P\to0$ in every $H^p$ by interpolation between
   $H^2$ and $H^\infty$ (and by inclusion when $p\le2$).
   Density of polynomials and uniform boundedness of $(T_N)$
   finish the proof.

   For arbitrary rotations, Lemma~\ref{lem:roots-of-unity} inserts
   unimodular factors $e^{-inj\theta_n}$ into the terms with
   exponent $m+nj$, leaving the orthogonality estimates unchanged.
   The maximal-function bounds are rotation invariant as well.
 \end{proof}


 \begin{theorem}[General Fej\'er phenomenon]   \label{thm:Fejer-on-subsequences}
   Suppose that $\Lambda\subset\DD$ is non-tangentially dense
   almost everywhere. Let $Q_n$ be the finite sampling operators
   constructed in \eqref{eq:Qn}. There exists a strictly increasing
   sequence $\vphi:\NN\to\NN$, depending only on the $Q_n$, such
   that
   \[
     \frac1N\sum_{j=1}^NQ_{\vphi(j)}f\longrightarrow f     \qquad\text{in }H^p
   \]
   for every $f\in H^p$ and every $1<p<\infty$.
   The indices can be selected inductively by finitely many
   $H^2$ scalar-product tests at each step. The associated Ces\`aro
   coefficient vectors remain uniformly bounded in $\CS_\Lambda$:
   \[
     \sup_N\Bignorm{\frac1N\sum_{j=1}^N        c^{(\vphi(j))}(f)}_{\CS_\Lambda}        \lesssim_{p,A}\norm{f}_{H^p}.
   \]
 \end{theorem}
 \begin{proof}
   Set $u_n^{(k)}:=Q_nz^k-z^k$ for $n\ge2$ and $k\ge0$.
   By Proposition~\ref{prop:weak-convergence}, for every fixed $k$
   we have $u_n^{(k)}\rightharpoonup0$ in $H^2$, while
   \[
     \sup_{n\ge2,k\ge0}\norm{u_n^{(k)}}_{H^2}\le M
   \]
   for some $M<\infty$. Choose $\vphi(1)\ge2$ arbitrarily.
   Having chosen $\vphi(1),\ldots,\vphi(j-1)$, take
   $\vphi(j)$ to be the smallest integer $\ell>\vphi(j-1)$
   for which
   \begin{equation}\label{eq:dyadic-decay}
     \bigl|\langle u_{\ell}^{(k)},
                      u_{\vphi(i)}^{(k)}\rangle_{H^2}\bigr|
       < 2^{-j},\qquad 0\le k\le j,\quad 1\le i<j.
   \end{equation}
   Such an $\ell$ exists: there are only finitely many conditions,
   and each scalar product tends to zero as $\ell\to\infty$.

   \proofskip Fix $k$. For all $j\ge\max\{2,k\}$, the new vector
   $u_{\vphi(j)}^{(k)}$ satisfies \eqref{eq:dyadic-decay}
   against every earlier selected vector. The finitely many
   remaining pairs contribute a constant depending on $k$, while
   \(\sum_{j\ge1}(j-1)2^{-j}<\infty\). Thus
   \[
     \Bignorm{\frac1N\sum_{j=1}^Nu_{\vphi(j)}^{(k)}}_{H^2}^2
       \le\frac{M^2}{N}+\frac{C_k}{N^2}\longrightarrow0.
   \]
   By linearity, the Ces\`aro averages converge strongly on
   polynomials in $H^2$.

   \proofskip    To pass to any $H^p$, fix $p>2$ and choose $s>p$.
   Proposition~\ref{prop:weak-convergence} gives uniform
   boundedness of $Q_n$ on $H^s$, hence of their Ces\`aro
   averages. Interpolation between strong convergence in $H^2$
   and boundedness in $H^s$ yields strong convergence for each
   polynomial in $H^p$. The case $1<p\le2$ follows directly from
   the $H^2$ result. Uniform boundedness on $H^p$ and density of
   polynomials finish the argument. The selected subsequence is
   independent of both $f$ and $p$.
 \end{proof}

 The conditions \eqref{eq:dyadic-decay} are finitely many,
 computable tests at each stage when the points of $\Lambda$ and the quadrature
 weights are given effectively. Recall indeed,
 \[
   Q_n(f)(z) = \sum_{j=1}^{J_n} m(E_{n,j}) \, f(\lambda_{n,j})  \; k_{\lambda_{n,j}}(z)   
 \]
 and the Hilbert-space inner products are computable from the kernel
 Gram matrix, for instance
 $[ k_\lambda,k_\mu ]_{H^2} = (1-\overline\lambda\,\mu)^{-1}$.  So at
 every step, the tests involve only finite sums of known
 quantities. We do not need to know how rapidly weak convergence takes
 place: it guarantees that the search eventually
 terminates. Theorem~\ref{thm:Fejer-on-subsequences} is not an
 unconditional computability theorem though, since for arbitrary
 sequences $\Lambda$ the initial construction of the finite shadow
 coverings need not be explicit; the proof of
 Theorem~\ref{thm:Fejer-on-subsequences} is therefore only algorithmic
 relative to an effective enumeration of $\Lambda$ and computable
 quadrature data.

 \subsection*{Acknowledgements and use of artificial intelligence. }

 This work was started in 2024 while visiting Markus Haase. We a
 grateful for finanacial help by the Christian-Albrechts-Universität
 zu  Kiel. \\[1ex]
 After an initial draft had been written, the author used OpenAI's
 ChatGPT to critically discuss proof strategies, and improve the
 presentation. The developing mathematical dialogue contributed
 substantially to the present version of the manuscript.  The author
 takes full responsibility for the mathematical content and
 correctness of the paper.

 \appendix

\section{Two auxiliary kernel calculations}

The first identity is the roots-of-unity filter in a form that covers
the examples and arbitrary rotations of radial blocks.

\begin{lemma}\label{lem:roots-of-unity}
  Let $n\ge1$, $0\le r<1$, $\theta\in\RR$, and
  $\omega=e^{2\pi i/n}$. Define
  \[
    Q_{n,r,\theta}f(z)      :=\frac1n\sum_{j=0}^{n-1}       f(re^{i\theta}\omega^j)\,k_{re^{i\theta}\omega^j}(z).
  \]
  For $m\ge0$ write $m=bn+s$, $0\le s<n$. Then
  \begin{equation}\label{eq:roots-of-unity-filter}
    Q_{n,r,\theta}(z^m)(z)      =\frac{r^{m+s}e^{i(m-s)\theta}z^s}            {1-r^ne^{-in\theta}z^n}.
  \end{equation}
  In particular, for $n>m$ the numerator is $r^{2m}z^m$,
  while $Q_{n,r,\theta} \eins = \frac{1}{1-r^ne^{-in\theta}z^n}$.
\end{lemma}
\begin{proof}
  Expanding each kernel in its geometric series, the coefficient
  of $z^\ell$ in $Q_{n,r,\theta}(z^m)$ is
  \[
    r^{m+\ell}e^{i(m-\ell)\theta} \; \Bigl(      \frac1n\sum_{j=0}^{n-1}\omega^{j(m-\ell)} \Bigr).
  \]
  The last factor is $1$ if $\ell\equiv m\pmod n$ and $0$
  otherwise. Summing over $\ell=s+jn$, $j\ge0$, gives
  \eqref{eq:roots-of-unity-filter}.
\end{proof}

The following standard estimate is recorded for convenience.
\begin{lemma}\label{lem:Hp-norms}
  Let $1<p<\infty$, $1/p+1/q=1$, and $N\ge1$.
  Uniformly in $N$, as $|\lambda|\uparrow1$,
  \[
    \Bignorm{\frac1{1-\overline\lambda \,z^N}}_{H^p}
       \asymp_p (1-|\lambda|)^{-1/q}.
  \]
\end{lemma}
\begin{proof}
  Write $r=|\lambda|$ and assume $r>1/2$.
  Periodicity and a change of variables show that
  \[
    \int_{\TT}|1-\overline\lambda\,\zeta^N|^{-p}\,dm(\zeta)
      =\frac1{2\pi}\int_{-\pi}^{\pi}
         \bigl((1-r)^2+2r(1-\cos t)\bigr)^{-p/2}\,dt,
  \]
  independently of $N$. Since $1-\cos t\asymp t^2$ for
  $|t|\le\pi$ and $p>1$, the last integral is comparable to
  \[
    \int_0^\pi\bigl((1-r)^2+t^2\bigr)^{-p/2}\,dt
       \asymp_p (1-r)^{1-p}.
  \]
  Taking $p$-th roots proves the assertion.
\end{proof}

\def\SUBMITTED{Submitted}
\def\TOAPPEAR{To appear in }
\def\PREPARATION{In preparation }

\def\cprime{$'$} \def\cprime{$'$}
\providecommand{\bysame}{\leavevmode\hbox to3em{\hrulefill}\thinspace}

\end{document}